\documentclass[12pt]{article}
\usepackage{mypaper}
\usepackage{amssymb}
\usepackage{amsmath}
\usepackage{amsthm}
\usepackage{epsfig}
\usepackage{hyperref}
\newtheorem{theorem}{Theorem}

\newenvironment{packed_enumerate}{
\setlength{\parsep}{0pt}
\setlength{\parskip}{0pt}
\begin{enumerate}
  \setlength{\itemsep}{1pt}
  \setlength{\parsep}{0pt}
  \setlength{\parskip}{0pt}
}{\end{enumerate}}

\title{The mathematics of Septoku}

\author{George I. Bell\\
\small \texttt{gibell@comcast.net}, \url{http://www.gibell.net/}}

\date{\dateline{Jan 21, 2008}\\
\small Mathematics Subject Classifications: 00A08, 97A20}

\begin{document} 
\maketitle

\begin{abstract}
Septoku is a Sudoku variant invented by Bruce Oberg,
played on a hexagonal grid of 37 cells.
We show that up to rotations, reflections, and
symbol permutations, there are only six valid
Septoku boards.
In order to have a unique solution, we show that
the minimum number of given values is six.
We generalize the puzzle to other board shapes,
and devise a puzzle on a star-shaped
board with 73 cells with six givens which has a unique solution.
We show how this puzzle relates to the unsolved
Hadwiger-Nelson problem in combinatorial geometry.
\end{abstract}

Sudoku
is one of the most popular pencil and paper puzzles,
and has sprouted a seemingly endless series of variations.
Witness the many recent books chock-full of exciting new variants
on the Sudoku theme \cite{Tridoku,Battle,Mutant}.

One variation is to consider the puzzle on a hexagonal rather
than square grid.
Bruce Oberg invented such a variation in 2006 and named it ``Septoku" \cite{ObergPuz}.
The name derives from the fact that seven is a significant number for this puzzle.
Septoku fit nicely into the seven theme at the 2006 conference ``Gathering 4 Gardner 7",
where it was Bruce Oberg's exchange gift \cite{ObergPuz}.

\begin{figure}[htb]
\centering
\epsfig{file=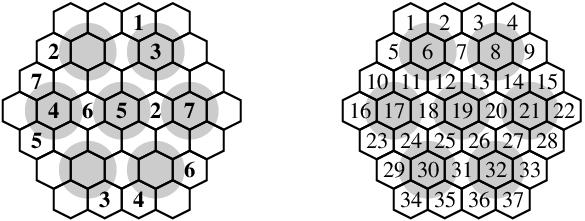}
\caption{A Septoku puzzle \cite{ObergPuz} and the Septoku board with cells numbered.}
\label{fig1}
\end{figure}

Figure~\ref{fig1}a shows a sample Septoku puzzle,
and Figure~\ref{fig1}b the numbering of the 37 hexagonal \textbf{cells}.
We will refer to a \textbf{row} of the board to include not only
the horizontal rows, but also the ``diagonal rows" parallel
to the other sides of the board.
For example the cells numbered $\{1$,$2$,$3$,$4\}$ form one row,
as do the cells $\{1$,$5$,$10$,$16\}$ and $\{3$,$8$,$14$,$21$,$28\}$.
We will abbreviate these as
``row 1$\rightarrow$4", ``row 1$\rightarrow$16" and ``row 3$\rightarrow$28".
In general, any subset of the board
that must contain distinct symbols will be called
a \textbf{region}.

The puzzle is to fill the 37 cells with
\textbf{seven} symbols (the numbers 1 to 7) such that
\begin{packed_enumerate}
\item Each row contains \textit{different} symbols.
There are seven rows in each of three directions, so 21 rows in all.
Note that most rows contain \textit{less} than 7 cells.
\item The seven \textbf{circular regions}
(hexagons of side two) must each contain all seven symbols.
For example the cells $\{1,$$2$,$5$,$6$,$7$,$11$,$12\}$ must
contain all of the numbers 1 to 7.
These circular regions are denoted by
gray circles as a visual aid to the
solver\footnote{This helpful suggestion is due to Nick Baxter.}.
Note that the circles overlap.
\end{packed_enumerate}
Henceforth, we will refer to a circular region as a \textbf{circle};
a particular circle is identified by its center,
i.e. ``circle 19" refers to cells
$\{12$,$13$,$18$,$19$,$20$,$25$,$26\}$.

A puzzle is defined by a blank board with some cells filled in,
called \textbf{givens} (Figure~\ref{fig1}a).
The puzzle is to fill out the remaining cells so that the above
rules are satisfied,
such a completed puzzle will be called a \textbf{valid} board.
A well-crafted puzzle has only one solution.

\begin{figure}[htb]
\centering
\epsfig{file=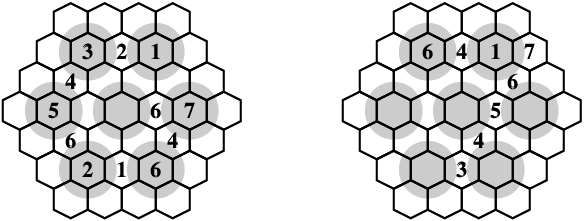}
\caption{``Medium" and ``hard" Septoku puzzles.}
\label{fig2}
\end{figure}

Before reading further, the reader should try the Septoku puzzle in Figure~\ref{fig1},
and the two harder puzzles in Figure~\ref{fig2}.
To get started with the Figure~\ref{fig1} puzzle, one should be able to deduce
the value in cells 13, 24, and 28, and then continue with cells 12, 11, 16, 25, \ldots \
For a collection of 28 Septoku puzzles that can be downloaded and printed, see \cite{ObergPuz}.

\section*{Seven Septoku Theorems}

On the face of it, solving Septoku seems similar to Sudoku.
But first impressions can be deceiving.
In fact, Septoku differs significantly
because it is much more tightly constrained.
As we will see, Septoku is also much easier to analyze.

In Sudoku, each cell (usually called a square in this case)
is in exactly three regions, a (normal) row, column, and $3\times 3$ block.
There are nine rows, nine columns, and nine $3\times 3$ blocks,
so 27 regions in all.
In Septoku there are 28 regions
(seven rows in each of three directions,
plus seven circles),
but the board is less than half the size of a Sudoku board.
Each cell is in at least four regions and twelve cells are in five regions.

When Bruce Oberg created the first Sudoku puzzles,
he noticed that all the valid boards he generated satisfied additional constraints
that were not required by the puzzle rules.
Here we see that these extra constraints on valid boards are easy to prove.

\begin{theorem}
In a valid Septoku board, each circle center must have a unique symbol.
In other words, $R_1=\{6$,$8$,$17$,$19$,$21$,$30$,$32\}$ must be a region
(contain each symbol 1--7).
\label{th:cc}
\end{theorem}

\begin{proof}
Suppose otherwise.
Then there must be some symbol that occurs twice
among the circle centers.
The center symbol is clearly unique by the rules
of the puzzle, so without loss of generality,
assume a 1 appears in cells 6 and 21.
There must be a 1 in the center circle 19,
and this can only be in cell 13 or 25.
But if it is in 25, then circle 8 can contain
no 1 at all.
Thus there must be a 1 in cell 13.
There must be a 1 in circle 17,
and this can now only be in cell 24,
and a 1 in circle 32, which can only be in cell 31.
But now circle 30 contains two 1's,
so this is not a valid Septoku board.
\end{proof}

\begin{theorem}
In a valid Septoku board, each corner and the center cell $19$ must have a unique symbol.
In other words, $R_2=\{1$,$4$,$16$,$19$,$22$,$34$,$37\}$ must be a region
(contain each symbol 1--7).
\label{th:corner}
\end{theorem}

\begin{proof}
Suppose otherwise.
Then there must be some symbol that occurs twice among the
corners (the center symbol is clearly unique by the
rules of the puzzle).
Without loss of generality, assume a 1 appears in cells 1 and 22.
By Theorem~\ref{th:cc} there must be a 1 in some circle center,
and this can only be in cell 8 or 30.
In either case, there can no longer be a 1
in the center circle 19, so this is not a valid Septoku board.
\end{proof}

\begin{theorem}
In a valid Septoku board,
each symbol must appear at least five times.
\label{th:5times}
\end{theorem}

\begin{proof}
Suppose the symbol 1 appears four or fewer times.
By Theorems~\ref{th:cc} and \ref{th:corner} each symbol must appear in a circle center,
and in a corner (or once in the center, 19).
If 1 occurs in a corner and one circle center,
then this covers only two of the circles out of seven.
There is no way the last two 1's can cover the remaining
five circles, since one cell is in common with at most two circles.
If 1 is in the center, then to cover the remaining 6 circles we must have a 1
in cells 7, 24 and 27 (or rotational equivalents),
but now we have two 1's in row 3$\rightarrow$29.
\end{proof}

If each symbol occurs 5 times this covers only 35 of the 37 cells,
so at least one symbol must occur more than 5 times.
In fact there are only two possibilities:
\begin{packed_enumerate}
\item Five of the symbols appear 5 times, two symbols appear 6 times.
\item Six of the symbols appear 5 times, one symbol appears 7 times.
\end{packed_enumerate}
We shall see that both of these possibilities can occur.

\begin{theorem}
Not counting reflections, rotations,
and permutations, there are only \textbf{six} valid Septoku boards.
Counting reflections, rotations, and permutations
as separate boards, there are $120,960 = 24\times 7!$ valid Septoku boards.
\label{th:ct}
\end{theorem}

\begin{figure}[htb]
\centering
\epsfig{file=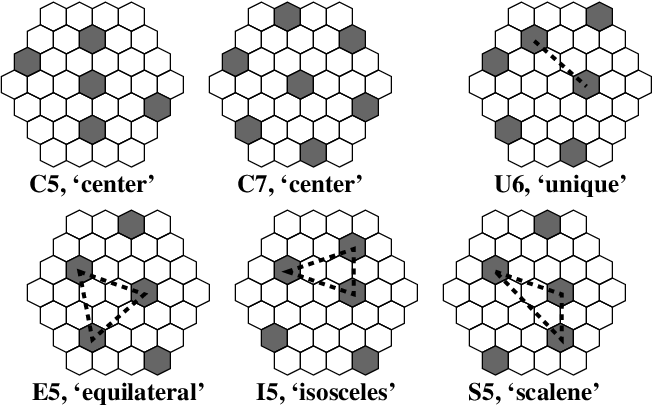}
\caption{Patterns for a number on a Septoku board:
the two patterns including the center, C5 and C7; the unique 6-cell pattern U6;
the 5-cell patterns E5, I5 and S5.}
\label{fig4}
\end{figure}

The central cell is the key to this board,
because whatever number is placed in cell 19,
there are only two patterns for this number over the rest of the board
(up to rotations and reflections).
These two patterns are shown in Figure~\ref{fig4},
where the dark-shaded cells identify all cells containing the same number.
To derive the two central patterns, start with the central cell 19, and suppose
the pattern includes at least one cell from the interior cells $I=\{7$,$11$,$14$,$24$,$27$,$31\}$.
Without loss of generality, we can assume this cell is 7.
We must cover circles 17 and 21, so the pattern must include cells 10 and 28,
or 15 and 23, and then the only way to cover circles 30 and 32 is by including cell 31.
This gives the pattern C5.
If we assume the pattern includes no cell in $I$ we get C7.
The two patterns C5 and C7 have $180^\circ$ and $60^\circ$ rotational symmetry,
respectively.
This is significant because (as we will see)
this symmetry is carried into the full solution.

Surprisingly, there is only one way to place six of the same numbers on
the board (up to symmetry),
and there are exactly three ways to place five of the same numbers on the board
(not including the center pattern C5).
These patterns are shown in  Figure~\ref{fig4}.
The simplest way to derive these patterns is from the center outward:
without loss of generality
assume the pattern includes cell 20.
We then need to have a number on circle 17,
and this can only be at cell 10 or 11
(or equivalently by reflection, cell 23 or 24).
Cell 10 leads to the unique 6-cell pattern (U6),
and cell 11 leads to the three 5-cell patterns: E5, I5, and S5.
A potential fourth pattern is obtained from I5 by moving a shaded cell from 8 to 4,
but this pattern violates Theorem~\ref{th:corner}, so is eliminated.
These patterns may be easily identified by the shapes formed by
their interior cells: ``unique" U6 (two interior cells), ``equilateral" E5, ``isosceles" I5,
and ``scalene" S5.
These shapes are useful for identifying these patterns in Septoku solutions.

Just as polyominoes are obtained by joining squares, a polyhex is obtained
by joining hexagons.
By looking at the patterns formed by all occurrences of a number,
we can view a Septoku problem as a polyhex packing problem!
Our polyhex pieces are the patterns shown in Figure~\ref{fig4}---these are
unusual polyhex pieces in that they are \textit{totally disconnected}.
Nonetheless, finding a valid Septoku board is equivalent to packing seven
of these polyhexes into the board.
We can rotate or flip any of the polyhexes, but we are not allowed to translate them.

\begin{figure}[htb]
\centering
\epsfig{file=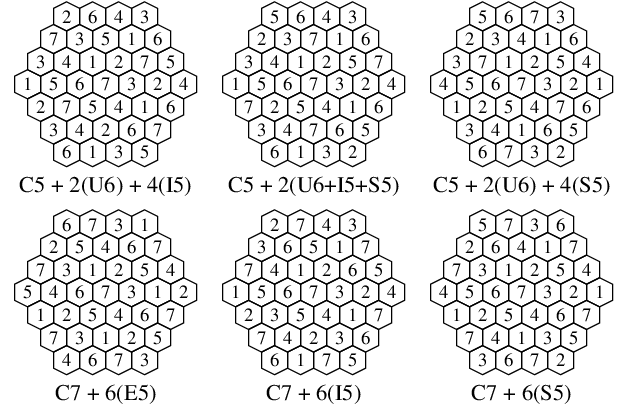}
\caption{The six solution boards, with the seven polyhex patterns each
is composed of given below.}
\label{fig5}
\end{figure}

Consider the case where two numbers appear six times,
and all other numbers five times.
We take two copies of the polyhex U6, the central polyhex C5,
and up to four copies of any of E5, I5 or S5.
Suppose we take U6 together with a reflection and rotation of U6---it is either
not possible to place the second U6, or not possible to place C5.
The only possibility is a copy of U6 together with a rotation of U6 by $180^\circ$.
The pattern C5 can then be fit in any of three orientations,
and each leads to a unique solution, given in the first row of Figure~\ref{fig5}.

The other case is to use the polyhex C7 together with six
of E5, I5 and/or S5.
It is not hard to see that these can only work out using C7 plus six copies
of the same polyhex, leading to another three unique solutions given in the
second row of Figure~\ref{fig5}.

The six Septoku solutions in Figure~\ref{fig5} satisfy some interesting
symmetry relations.
We define $(Rot)$ as the transformation which rotates the board
clockwise $60^\circ$.  We use standard notation to represent a
permutation as a product of cycles.
The cycle $(14)$ describes the permutation which sends
$1\mapsto 4$, $4\mapsto 1$, and leaves all other numbers alone.

We can easily verify that all solution boards $B$ in Figure~\ref{fig5} satisfy
\begin{equation}
B(Rot)^3(14)(25)(36) = B.
\label{eq:180deg}
\end{equation}
The boards in the second row of Figure~\ref{fig5} satisfy an even stricter relation
\begin{equation}
B(Rot)(123456) = B.
\label{eq:60deg}
\end{equation}

The lower-left board in Figure~\ref{fig5} has an additional symmetry:
all seven numbers are present in \textit{every possible circle on the board}!
This includes the seven Septoku circles plus twelve additional circles centered
at all remaining interior cells.
As we shall see, this ``all equilateral" pattern can be extended periodically to infinity.

If we take the six Septoku boards, and apply any of 12 symmetry transformations followed
by any of $7!$ permutations, we will generate $6(12)(7!)$ valid Septoku boards.
However, not all of these boards are distinct, due to the symmetry Equations
(\ref{eq:180deg}) and (\ref{eq:60deg}).
For any board in the top row of Figure~\ref{fig5}, the total number of boards is reduced by
a factor of 2, and in the bottom row by a factor of 6.
This reduces the total number of distinct boards to $3(12/2+12/6)(7!) = 24\times 7!$.

\begin{figure}[htb]
\centering
\epsfig{file=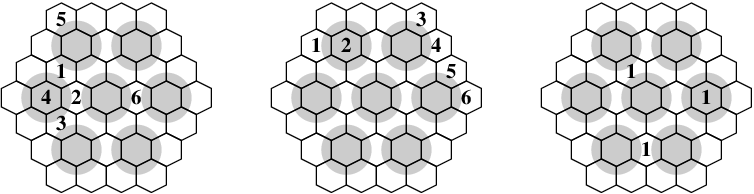}
\caption{Two puzzles with 6 givens and unique solutions,
and a puzzle with 3 givens and a unique solution up to a permutation.}
\label{fig6}
\end{figure}

\begin{theorem}
For a Septoku puzzle with a unique solution,
the minimum number of givens is six.
\label{th:givenmin}
\end{theorem}
\begin{proof}
Clearly the answer cannot be less than six.
If it was, there must be two of the seven numbers
not even among the givens.
If there is one solution, a different one can be obtained
by swapping these two numbers.

To show the answer is six, we need only find a puzzle
with six givens and a unique solution.
We wrote a program that checks puzzles for a unique
solution by matching against all possible symmetry
transformations and permutations of the six possible boards.
We found hundreds of six given puzzles with
unique solutions, two are shown in Figure~\ref{fig6}.
\end{proof}

To narrow the solution down to one of the six types shown
in Figure~\ref{fig5}, only 3 givens are required.
For example, note that the `equilateral' pattern E5 appears in only
one solution.  Thus, if we use the 3 givens on the right board in
Figure~\ref{fig6}, the solution can only be the ``all equilateral"
board on the lower left of Figure~\ref{fig5},
or any permutation of it which fixes 1.

The symmetry relations in Equation (\ref{eq:180deg})
generate a useful algorithm for solving Septoku boards.

\begin{theorem}
In a valid Septoku board,
the seven symbols can be partitioned into three pairs
and one unpaired symbol which is the value of the central cell 19.
Let $x$ and $y$ be any two cells that map to each other under
$180^\circ$ rotation of the board.
Then either
\begin{packed_enumerate}
\item $x$ and $y$ have different values taken from one of the three pairs, or
\item $x$ and $y$ both have the value of the central cell
(clearly this is only possible if $x$ and $y$ are not in the same row).
\end{packed_enumerate}
\label{th:big}
\end{theorem}

\begin{proof}
For the six solutions in Figure~\ref{fig5}, Theorem~\ref{th:big} is
just a restatement of Equation (\ref{eq:180deg}).
Since any solution is a symmetry transformation plus permutation
of one of these, all solutions satisfy Theorem~\ref{th:big}.
\end{proof}

Once we determine the pairings, Theorem~\ref{th:big} is very useful in
solving Septoku puzzles---because it says that once we
determine a number in a certain cell, we can immediately fill out the
cell diametrically opposite (mapped by $180^\circ$ rotation).

With these additional insights,
the reader may wish to revisit the puzzles in Figure~\ref{fig2}.
For example, in the ``hard" puzzle, one Theorem~\ref{th:big} pair is $\{3,4\}$,
and this implies that cell 12 must contain a 3.
The reader should find that completing this puzzle is not difficult using
Theorem~\ref{th:big}.
One puzzle that remains difficult is that in Figure~\ref{fig6}b.
In this puzzle no pairings can be determined from the starting givens.

The last theorem shows that the requirements on the central circle 19 can be removed,
and the puzzle is unchanged.

\begin{theorem}
Consider Septoku, but drop the requirement
that the symbols in the center circle 19 be distinct.
There are no solutions to the puzzle where the cells
in the center circle are not distinct.
\label{th:center}
\end{theorem}
\begin{proof}
Suppose a solution exists with two symbols
the same in circle 19.
Because the rows must be distinct,
the central number cannot be the duplicate, so
without loss of generality assume that the symbol 1 appears
in cells 13 and 26.
But circles 17 and 21 must both contain 1's,
and these can only be on row 16$\rightarrow$22.
There cannot be two 1's along this row,
so no solution is possible.
\end{proof}

\section*{Larger Boards}

\begin{figure}[htb]
\centering
\epsfig{file=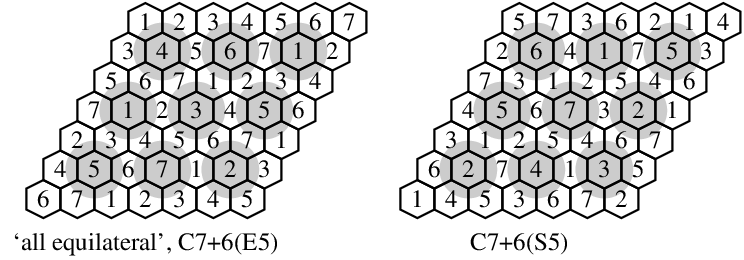}
\caption{The two valid rhombus Septoku boards.}
\label{fig10}
\end{figure}

When analyzing Septoku, I noticed that some Septoku solutions
have the property that the solution can be extended off the board,
while holding to a natural extension of the rules.
We can extend two opposite edges of the board to
create a 49-cell \textbf{rhombus} Septoku board in Figure~\ref{fig10}.

The rhombus (or diamond)
board has two sets of rows which all have 7 cells,
and another set of rows which have 2 to 7 cells.
On this board, every symbol must appear exactly seven times.
All the Figure~\ref{fig4} patterns can be extended to this board
\textit{except} for $U6$ and $I5$, and of the six solutions in Figure~\ref{fig5}
only two do not involve $U6$ or $I5$,
therefore there are only \textbf{two} valid rhombus Septoku boards (Figure~\ref{fig10}).

The ``all equilateral" solution in Figure~\ref{fig10} has additional symmetries in
that \textit{every circle} contains all seven numbers.
Extended periodically, this solution shows that \textit{any} Septoku board,
no matter how large or how the circles are arranged, has at least one valid board.

If we interpret each number as a color, the ``all equilateral" solution gives a
7-coloring of the plane which is well-known to graph theorists.
If the diameter of each hexagon is just under one, this coloring has the
remarkable property that any two points in the plane a distance exactly
one apart have different colors.
This coloring provides an upper bound for the so-called Hadwiger-Nelson problem,
which asks for the minimum number of colors for which
such a coloring of the plane is possible.
As shown in 1961, the minimum number of colors required lies
between 4 and 7 (inclusive) \cite{GardnerWL, Ogilvy}.
In 2018, the lower bound was increased to 5 \cite{deGrey, Exoo},
but the Hadwiger-Nelson problem remains a significant
unsolved problem in combinatorial geometry \cite{TOPP}.

\begin{figure}[htb]
\centering
\epsfig{file=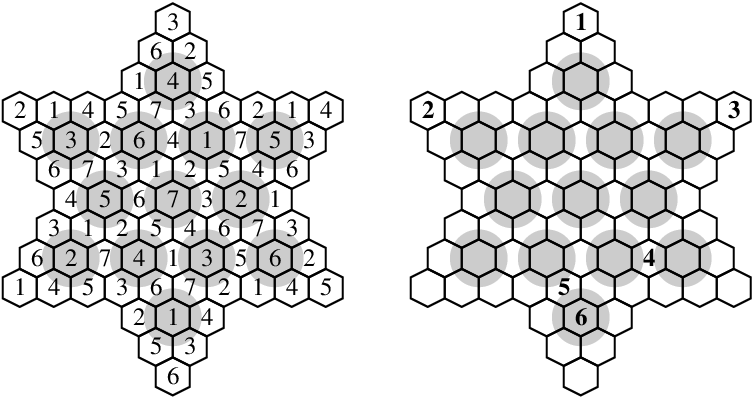}
\caption{A valid star Septoku board, a puzzle with 6 givens and a unique solution.}
\label{fig11}
\end{figure}
Another variation which is possible is the 73-cell
\textbf{star} Septoku board (Figure~\ref{fig11}).
On this board, some of the rows have more than seven cells.
We therefore modify the rules to specify that the symbols in
any connected subset of a row with seven or fewer cells must
contain distinct symbols.
There are only two valid star Septoku boards.

Not surprisingly, it is easy to create many
puzzles in rhombus or star Septoku with 6 givens that have unique solutions.
An example of one such puzzle is shown in Figure~\ref{fig11}.

\begin{figure}[htb]
\centering
\epsfig{file=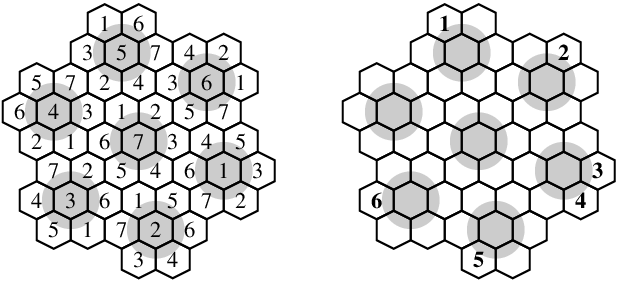}
\caption{A valid flower Septoku board, a puzzle with 6 givens and a unique solution.}
\label{fig12}
\end{figure}

In all the board types discussed so far, the circular regions
overlap.
It is possible to create a board on a hexagonal grid where the
circular regions do not overlap, more reminiscent of
ordinary Sudoku.
Figure~\ref{fig12} shows a \textbf{flower} Septoku board,
which has 49 cells and 34 regions.


Since the circular regions no longer overlap,
one might guess that flower Septoku is quite different from
normal Septoku.
However, a similar analysis reveals that there are only \textbf{three}
different flower Septoku boards (up to a symmetry transformation
and permutation).
These three solutions also satisfy Theorem~\ref{th:big}---so despite
the change in geometry,
solutions to flower Septoku have the same symmetry
properties as in regular Septoku.

\section*{Conclusions}

Although Septoku appears similar to Sudoku, we have found that
it has only six fundamentally different solution boards.
We have shown that a Septoku puzzle must have at least
six givens to have a unique solution, 
and have devised puzzles with six givens
that can be quite difficult to solve by hand.
Theorems~\ref{th:cc}, \ref{th:corner}, and especially \ref{th:big}
are useful for solving Septoku puzzles.

The main difference between Sudoku and Septoku are the
extra constraints imposed by an additional set of
parallel rows in a hexagonal grid.
In any variation of Septoku we have considered,
these additional constraints have resulted in
a relatively small number of solutions.


How might we modify the rules of Septoku so that there are more
than six solution boards?
One way is to reduce the number of regions,
but we have shown that removing the central region
does not change the puzzle.
Another option is to allow 8 symbols and keep all the regions the same.
If we do this, all our theorems are no longer valid,
and many more solution boards are possible.
However, it is now quite difficult to specify a unique
solution to a puzzle using a reasonable number of givens.
For example,
it is not hard to find an 8 symbol Septoku puzzle with 36 givens and
only one open cell which does not have a unique solution!

\section*{Acknowledgments}
I thank Bruce Oberg for inventing this fascinating puzzle, and
Ed Pegg for pointing out the connection to the Hadwiger-Nelson problem.
Thanks to Wei-Hwa Huang for sharing his analysis of Septoku,
and his idea that the two central patterns were the key
to the analysis.

\section*{Solutions}

\begin{figure}[ht]
\centering
\epsfig{file=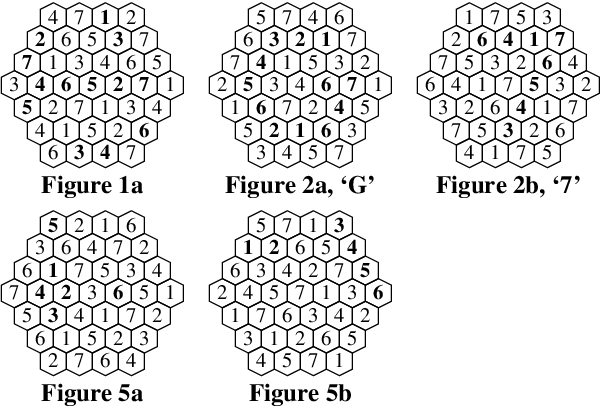}
\caption{Unique solutions for problems in figures as indicated.}
\label{fig9}
\end{figure}

See Figure~\ref{fig9} for solutions to the Septoku puzzles in
Figures~\ref{fig1}, \ref{fig2} and \ref{fig6}.
The unique solution to the star Septoku problem in Figure~\ref{fig11}b can be
found by applying the permutation $(16543)$ to the solution in Figure~\ref{fig11}a.
For the flower Septoku problem in Figure~\ref{fig12}b,
simply label each circular region with the same pattern,
that in the center circle of the board in Figure~\ref{fig12}a.


\end{document}